\documentstyle[12pt,amssymb]{amsart}

\title[Actions of Quantum Groups on Operator Algebras]
{\bf Ergodic Actions of Universal Quantum Groups on Operator Algebras}
\author[Shuzhou Wang]
{\bf Shuzhou Wang}
\address{Department of Mathematics, University of California,
Berkeley, CA 94720
\newline \indent
Fax: 510-642-8204
}
\email{szwang@@math.berkeley.edu}

\date{}

\newtheorem{DF}{Definition}[section]
\newtheorem{LM}[DF]{Lemma}
\newtheorem{PROP}[DF]{Proposition}
\newtheorem{TH}[DF]{Theorem}
\newtheorem{COR}[DF]{Corollary}
\newtheorem{RMK}[DF]{Remark}
\newtheorem{RMKS}[DF]{Remarks}
\newtheorem{PROB}[DF]{Problem}

\newcommand{\bgdf}{\begin{DF}}
\newcommand{\nddf}{\end{DF}}
\newcommand{\bglm}{\begin{LM}}
\newcommand{\ndlm}{\end{LM}}
\newcommand{\bgprop}{\begin{PROP}}
\newcommand{\ndprop}{\end{PROP}}
\newcommand{\bgth}{\begin{TH}}
\newcommand{\ndth}{\end{TH}}
\newcommand{\bgcor}{\begin{COR}}
\newcommand{\ndcor}{\end{COR}}
\newcommand{\bgrmk}{\begin{RMK}}
\newcommand{\ndrmk}{\end{RMK}}
\newcommand{\bgrmks}{\begin{RMKS}}
\newcommand{\ndrmks}{\end{RMKS}}
\newcommand{\bgprob}{\begin{PROB}}
\newcommand{\ndprob}{\end{PROB}}
\newcommand{\bgeq}{\begin{eqnarray}}
\newcommand{\ndeq}{\end{eqnarray}}
\newcommand{\bgeqq}{\begin{eqnarray*}}
\newcommand{\ndeqq}{\end{eqnarray*}}

\newcommand{\QED}{\hfill Q.E.D.} 
\newcommand{\vv}{\vspace{4mm}\\}

\numberwithin{equation}{section}

\newcommand{\dfref}[1]{Definition~\ref{#1}}

\newcommand{\propref}[1]{Proposition~\ref{#1}}
\newcommand{\thref}[1]{Theorem~\ref{#1}}

\newcommand{\corref}[1]{Corollary~\ref{#1}}

\newcommand{\secref}[1]{\S\ref{#1}}

\begin{document}

\begin{abstract}
We construct ergodic actions of compact quantum
groups on $C^*$-algebras and von Neumann algebras, and exhibit
phenomena of such actions that are of different nature from
ergodic actions of compact groups.

In particular, we construct:
(1). an ergodic action of the compact quantum
$A_u(Q)$ on the type III$_\lambda$ Powers factor $R_\lambda$
for an appropriate positive $Q \in GL(2, {\Bbb R})$;
(2). an ergodic action of the compact quantum group
$A_u(n)$ on the hyperfinite $\mathrm{II}_1$ factor $R$;
(3). an ergodic action of the compact quantum group $A_u(Q)$ on the
Cuntz algebra ${\cal O}_n$ for each positive matrix $Q \in GL(n, {\Bbb C})$;
(4). ergodic actions of compact quantum groups on the
their homogeneous spaces, as well as an example of a non-homogeneous
classical space that admits an ergodic action of a compact quantum group.
%
\end{abstract}

\maketitle

\section{
Introduction}

%
%

It is well known that compact groups
admit no ergodic actions on operator algebras other than the
finite ones (i.e. those with finite traces) \cite{HLS}.
Therefore, there arosed the following basic problem
(cf p76 of \cite{HLS}):
Construct an ergodic action of a semisimple compact Lie group on the
Murray-von Neumann $\mathrm{II}_1$ factor $R$. Later,
%
%
Wassermann developed some general theory of ergodic actions of
compact groups on operator algebras and
%
%
showed that $SU(2)$ cannot act ergodically on $R$
\cite{AWassermann1,AWassermann3},
leaving experts the doubt that semisimple compact Lie groups
admit ergodic actions on $R$ at all.
In \cite{Boca1}, Boca studied the general theory of
ergodic action of compact quantum groups \cite{Wor5} on $C^*$-algebras
and generalized some basic results on ergodic actions of
compact groups to compact quantum groups. But so far
there is still a lack of non-trivial examples of ergodic actions
of compact quantum groups on operator algebras.

The purpose of the present paper is two-fold, which is
in some sense opposite to that of Boca \cite{Boca1}. First, we show
that some new phenomena can occur for ergodic actions of
quantum groups. Second, we supply some
general methods to construct ergodic actions of compact
quantum groups on operator algebras and give several non-trivial examples of
such actions. We show that the universal compact matrix quantum
groups $A_u(Q)$ of \cite{W5,W1} admit ergodic actions on both
the (infinite) injective
factors of type III (for $Q \neq c I_n$, $c \in {\Bbb C}^*$)
and the (infinite) Cuntz algebras (for $Q>0$). We construct
an ergodic action of the universal compact matrix quantum group of Kac type
$A_u(n)$ on the hyperfinite factor $R$, which may not admit
ergodic actions of any semisimple compact Lie group \cite{AWassermann3}.
We also study ergodic actions of compact quantum groups
on their homogeneous spaces and show that there are non-homogeneous
{\em classical spaces} that admit ergodic actions of quantum groups.
These results show that compact quantum groups have a much
richer theory of ergodic actions on operator algebras than
compact (Lie) groups.

Unlike Boca \cite{Boca1}, we study actions of compact quantum groups
on both $C^*$-algebras and von Neumann algebras,
not just $C^*$-algebras. Our construction of
ergodic actions of compact quantum groups on von Neumann
algebras come from their ``measure preserving'' actions on
$C^*$-algebras, just as in the classical situation
(see \thref{inducedaction}). One of our constructions of ergodic actions
(see Sect. 3) uses tensor products of irreducible representations of compact
quantum groups. This method was first used by Wassermann \cite{AWassermann4}
in the setting of Lie groups (instead of quantum groups)
to construct subfactors from their ``product type actions''.
At the other extreme, actions of quantum groups with large fixed
point algebras (i.e. prime actions) have been studied by many authors,
see, e.g. \cite{Cuntz2,CDPR1}.
Generalizing the canonical action of compact Lie groups on the Cuntz
algebras \cite{Cuntz} introduced by Doplicher-Roberts \cite{DopRob6,Dop1},
Konishi et al \cite{KMW} study (the non-ergodic) action of
$SU_q(2)$ on the Cuntz algebra ${\cal O}_2$ and its CAR subalgebra and show
that their fixed point algebras coincide (see also \cite{Marciniak1}).
This result is extended to $SU_q(n)$ by Paolucci \cite{Paolucci1}.
In \cite{Nakagami1}, this action of the quantum group $SU_q(n)$ is
induced to a (non-ergodic) action on
the Powers factor $R_\lambda$ by a rather complicated method,
which follows from our result \thref{inducedaction} in a much simpler
and more conceptual manner.
%

The contents of this paper are as follows. In \secref{induce},
we give a general method of construction of quantum group actions
on von Neumann algebras from their ``measure preserving'' action
on $C^*$-algebras. Using this and a result of Banica \cite{Banica2}
on the tensor products of
the fundamental representation of $A_u(Q)$, we construct in
\secref{UHF} an ergodic action of a universal quantum groups $A_u(Q)$ on the
Powers factor $R_\lambda$ of type $\mathrm{III}_\lambda$ and
and an ergodic action of $A_u(n)$ on the hyperfinite
$\mathrm{II}_1$ factor $R$. In \secref{fixedpoint}, using results of
Banica \cite{Banica1}, we show that the
fixed point subalgebra of $R$ under the quantum subgroup $A_o(n)$
of $A_u(n)$ is also a factor and that the action of $A_o(n)$ on $R$ is prime.
In \secref{Cuntz},
we construct ergodic action of $A_u(Q)$ on the Cuntz algebras and on the
injective factor $R_{\infty}$ of type $\mathrm{III}_1$ as well as the other
factors of type $\mathrm{III}$. It is also shown that the (unimodular)
compact quantum group $A_u(n)$ of Kac type acts ergodically
on the injective factor of type $\mathrm{III}_{\frac{1}{n}}$,
a fact rather surprising to us.
In the last section \secref{quotient}, we
study ergodic actions of compact quantum groups on their ``quotient spaces'',
and show that the quantum automorphism group $A_{aut}(X_4)$ acts ergodically
on the classical space $X_4$ with $4$ points, but $X_4$ is not isomorphic
to a quotient space. We point out that instead of using the
fundamental representation of $A_u(Q)$, we can
also use representations of free products of compact quantum groups
\cite{W1} in the examples in \secref{UHF} and \secref{Cuntz} for the
constructions of ergodic actions.

\section{
\label{induce}
Lifting actions on $C^*$-algebras to von Neumann algebras}

In this section, we describe (\thref{inducedaction})
how to construct ergodic actions of compact quantum groups on
von Neumann algebras from ``measure preserving'' actions on
noncommutative topological spaces (i.e. $C^*$-algebras).

To fix notation, we first recall some basic notions concerning actions of
quantum groups on operator algebras (\cite{BS2,Boca1,Pod6,W15}).
For convenience in this paper, we will use the definition given in
\cite{W15} for the notion of actions of compact quantum groups on
$C^*$-algebras. As in \cite{W15},
Woronowicz Hopf $C^*$-algebras are assumed to be full in order to
define morphisms. We adapt the following convention
(see \cite{W1,W5,W15}): when $A=C(G)$ is a Woronowicz Hopf $C^*$-algebra,
we also say that $A$ is a compact quantum group,
referring to the dual object $G$.

\bgdf (cf \cite{W15})
\label{qact}
A (left) {\bf action} of a compact quantum group $A$ on a
$C^*$-algebra $B$ is a unital *-homomorphism  $\alpha$ from $B$ to
$B \otimes A$ such that

(1). $( id_B \otimes \Phi ) \alpha = ( \alpha \otimes id_A ) \alpha,$
where $\Phi$ is the coproduct on $A$;

(2). $(id_B \otimes \epsilon) \alpha = id_B$,
where $\epsilon$ is the counit on $A$;

(3). There is a dense *-subalgebra $\cal B$ of $B$, such that
$$\alpha ({\cal B}) \subseteq {\cal B} \otimes {\cal A},$$
where $\cal A$ is the canonical dense *-subalgebra of $A$.
\nddf
\noindent
{\em Remarks.}
(1).
The definition above is equivalent to the one in Podles \cite{Pod6}.
As in \cite{Pod6}, we do not impose the condition that $\alpha$ is
injective, which is required in \cite{BS2,Boca1},
though the examples constructed in this paper satisfy
this condition. We conjecture that this condition is
a consequence of the other conditions in the definition.
A special case of this conjecture says that
the coproduct of a Woronowicz Hopf $C^*$-algebra
is injective, which is true for both the full Woronowicz Hopf
$C^*$-algebras (because of the counital property) and
the reduced ones (because of Baaj-Skandalis \cite{BS2}).
Even if this conjecture is false, one can
still obtain an injective $\tilde{\alpha}$ from $\alpha$
by passing to the quotient of $B$ by the kernel of $\alpha$.
We leave the verification of the latter as an exercise for the reader.
%

(2).
The above notion of left action of quantum group $G$ would be
called right coaction of the Woronowicz Hopf $C^*$-algebra $C(G)$ by
some other authors. But we prefer the more geometric term ``action of
quantum group''.
We can similarly define a right action of a quantum group $G$, which
would be called a left coaction of the Woronowicz Hopf
$C^*$-algebra $C(G)$ by some other specialists.


\bgdf
Let $\alpha$ be an action of a compact quantum group $A$ on $B$.
An element $b$ of $B$ is said to be {\bf fixed under $\alpha$}
(or {\bf invariant under $\alpha$}) if
\bgeq
\alpha (b) = b \otimes 1_A.
\ndeq
The {\bf fixed point algebra} $B^\alpha$
(or $B^A$ if no confusion arises) of the action $\alpha$ is
\bgeq
B^\alpha = \{ b \in B \; | \; \alpha (b) = b \otimes 1_A \}.
\ndeq
The action of $A$ is said to be {\bf ergodic}
if $B^\alpha = {\Bbb C} I$.
A continuous functional $\phi$ on $B$ is said to be
{\bf invariant under $\alpha$} if
\bgeq
 (\phi \otimes id_A) \alpha (b) = \phi (b) 1_A.
\ndeq
\nddf

Fix an action $\alpha $ of a compact quantum group $A$ on $B$. Let $h$
be the Haar state on $A$ \cite{Wor5,W1,Daele4}. Then we have

\bgprop
\label{denseB}

(1). The map $E = (1 \otimes h) \alpha$ is a projection of norm one
from $B$ onto $B^\alpha$;

(2). Let
\bgeq
{\cal B}^\alpha = \{ b \in {\cal B} \; | \; \alpha (b) = b \otimes 1_A \}.
\ndeq
Then ${\cal B}^\alpha$ is norm dense in $B^\alpha$. Hence the
action $\alpha$ is ergodic if and only if it is so when restricted to
the dense *-subalgebra ${\cal B}$ of $B$.
\ndprop
\pf
(1).
This is an easy consequence of the following form of the invariance of the
Haar state (cf \cite{Wor5}):
$$(id_A \otimes h) \Phi (a) = h(a) 1_A, \; \; \; a \in A.$$

(2).
If $b \in B^\alpha$, then $b$ can be approximated in
norm by a sequence of elements $b_l \in \cal B$. Let ${\bar b}_l$ be the
average of $b_l$:
$$ {\bar b}_l = (1_B \otimes h) \alpha (b_l).$$
Then from part (1) of the proposition, $\bar{b}_l \in B^\alpha$. From
condition (3) of \dfref{qact}, we see that
$\bar{b}_l \in {\cal B}^\alpha$.
Moreover,
\bgeqq
\|{\bar b}_l - b \|
&=& \| (1 \otimes h) \alpha(b_l - b)\| \\
&\leq& \| (1 \otimes h) \alpha \| \|b_l - b\| \rightarrow 0.
\ndeqq
The rest is clear.
\QED
\vv
Preserve the notation above. Let $\frak A$ be the von Neumann algebra
generated by the GNS representation $\pi_h$ of $A$ for the state $h$.
Then $\frak A$ is a Hopf von Neumann algebra. For later use,
we need to adapt the definitions above to the situation of
von Neumann algebras. 

\bgdf
\label{vonNeumanncoact}
A {\bf right coaction} of a Hopf von Neumann algebra $\frak A$ on a
von Neumann algebra $\frak B$ is a normal homomorphism $\alpha$ from
$\frak B$ to ${\frak B} \otimes {\frak A}$ such that

(1). $ ( id_{\frak B} \otimes \Phi ) \alpha =
( \alpha \otimes id_{\frak A} ) \alpha, $
where $\Phi$ is the coproduct on $\frak A$;

(2). $ \alpha ({\frak B}) (1 \otimes {\frak A} ) $
generates the von Neumann algebra ${\frak B} \otimes {\frak A}$.
\nddf

The main reason why we use the term ``coaction of Hopf von Neumann algebra''
is that von Neumann algebras are measure-theoretic objects
instead of geometric-topological objects (cf. Remark (2) after \dfref{qact}).
Condition (2) in the above definition is an analog of the density
condition as used in the Hopf $C^*$-algebra setting \cite{BS2,Pod6}.
It is well known that there is no analogue of counit in
the Hopf von Neumann algebra situation simply because a
von Neumann algebra corresponds to
a measure space in the commutative case (the simplest case), and functions
are defined only up to sets of measure zero. Hence we do not have an
analog of condition (2) of \dfref{qact} for Hopf von Neumann algebras
coactions.

If $\frak A$ comes from the GNS-representation of the Haar state
on a compact quantum group $A$ and $\frak A$
coacts on the right on some von Neumann algebra $\frak B$,
we will abuse the terminology by saying that the quantum group
$A$ acts on $\frak B$.
Other notions such as invariant elements (or functionals), fixed point
algebra and ergodic actions in the
$C^*$-case above can also be carried over to to the
von Neumann algebra situation.

The main result of this section is the following

\bgth
\label{inducedaction}
Let $B$ be a $C^*$-algebra endowed with an action $\alpha$ of a compact
quantum group $A$. Let $\tau$ be an $\alpha$-invariant state on $B$.
Then

(1).
$\alpha$ lifts to a coaction $\tilde{\alpha}$ of the Hopf von Neumann
algebra ${\frak A} = \pi_h (A) ^{\prime \prime}$ on the von Neumann algebra
${\frak B} = \pi_{\tau}(B) ^{\prime \prime}$ defined by
\bgeq
\tilde{\alpha} (\pi_\tau (b)) = (\pi_\tau \otimes \pi_h) \alpha(b),
\hspace{1cm} b \in B,
\ndeq
where $\pi_h$ and $\pi_\tau$ are respectively the GNS representations
associated with the Haar state $h$ on $A$ and the state $\tau$ on $B$.

(2).
If $\alpha$ is ergodic, then so is $\tilde{\alpha}$.
\ndth
\pf
(1).
We will only show that the natural map $\tilde{\alpha}$ given
on the dense subalgebra $\pi_\tau(B)$ by
$$\tilde{\alpha} (\pi_\tau (b)) = (\pi_\tau \otimes \pi_h) \alpha(b),
\hspace{1cm} b \in B,$$
is well defined and extends to a {\em normal morphism} from $\frak B$
to ${\frak B} \otimes {\frak A}$.

Let $b \in B$ and $a \in A$. Denote by $\tilde b$ and
$\tilde a$ respectively the corresponding elements of the Hilbert spaces
$H=L^2(B, \tau)$ and $K = L^2(A, h)$. Define an operator $U$ on
$H \otimes K$ by
\bgeq
U(\tilde{b} \otimes \tilde{a}) = (\pi_\tau \otimes \pi_h) \alpha(b)
(\tilde{1}_B \otimes \tilde{a}).
\ndeq
Then since $\tau$ is $\alpha$ invariant, we have
\bgeqq
<U(\tilde{b} \otimes \tilde{a}), U(\tilde{b} \otimes \tilde{a})>
&=& (\tau \otimes h) (1_B \otimes a^*) \alpha(b^* b) (1_B \otimes a) \\
&=& aha^* ((\tau \otimes id_A) \alpha(b^* b))
= aha^* (\tau (b^* b) 1_A) \\
&=& <\tilde{b} \otimes \tilde{a}, \tilde{b} \otimes \tilde{a}>,
\ndeqq
where $aha^*$ is the functional on $A$ defined by
$$aha^*(x) = h(a^* x a), \hspace{1cm} x \in A.$$
Hence $U$ is an isometry. Since
$\alpha (B) (1 \otimes A)$ is dense in $B \otimes A,$
$U$ is a unitary operator. We also have
\bgeqq
(\pi_\tau \otimes \pi_h) \alpha (b) U (\tilde{b'} \otimes \tilde{a'})
&=& (\pi_\tau \otimes \pi_h) \alpha (b)
(\pi_\tau \otimes \pi_h) \alpha (b') (\tilde{1}_B \otimes \tilde{a'})  \\
&=& (\pi_\tau \otimes \pi_h) \alpha (b b') (\tilde{1}_B \otimes \tilde{a'})
= U (\pi_\tau(b) \tilde{b'} \otimes \tilde{a'}) \\
&=& U (\pi_\tau(b) \otimes 1) (\tilde{b'} \otimes \tilde{a'}).
\ndeqq
That is
\bgeq
\label{covariantaction1}
(\pi_\tau \otimes \pi_h) \alpha(b) = U (\pi_\tau(b) \otimes 1) U^*.
\ndeq
Condition (1) of \dfref{vonNeumanncoact} follows immediately.
Since $ \alpha (B) (1 \otimes A ) $
is dense in $B \otimes A$ (cf. Remark (1) after \dfref{qact}
and Podles \cite{Pod6}),
Condition (2) of \dfref{vonNeumanncoact} follows.

(2).
Assume $\alpha$ is ergodic. Let $z \in {\frak B}$ be a fixed element
under $\tilde{\alpha}$:
$$\tilde{\alpha} (z) = z \otimes 1_{\frak A}.$$
Let $b_n \in  B$ be a net of elements such that $\pi_\tau(b_n) \rightarrow z$
in the weak operator topology. Consider the average of $\pi_\tau(b_n)$
integrated over the quantum group $A$,
$$z_n = (id_{\frak B} \otimes h) \tilde{\alpha} (\pi_\tau (b_n)),$$
where we use the same letter $h$ to denote the faithful normal state
on $\frak A$ determined by the Haar state $h$ on $A$.
Then one can verify that $z_n \rightarrow z$ in the weak operator topology.
Moreover, using
$$(id_{\frak A} \otimes h) \Phi (a) = h(a) 1_{\frak A}, \; \; \; a \in
{\frak A},$$
where we denote the coproduct on $\frak A$ by the same symbol
as the coproduct $\Phi$ on $A$, we have
\bgeqq
\tilde{\alpha} (z_n)
&=& (id_{\frak B} \otimes id_{\frak A} \otimes h)
(\tilde{\alpha} \otimes id_{\frak A}) \tilde{\alpha}(\pi_\tau(b_n)) \\
&=& (id_{\frak B} \otimes id_{\frak A} \otimes h) (id_{\frak B} \otimes \Phi )
\tilde{\alpha} (\pi_\tau(b_n)) \\
&=& (id_{\frak B} \otimes (id_{\frak A} \otimes h) \Phi )
\tilde{\alpha}(\pi_\tau(b_n)) \\
&=& (id_{\frak B} \otimes h) \tilde{\alpha}(\pi_\tau(b_n)) \otimes 1_{\frak A}
= z_n \otimes 1_{\frak A}.
\ndeqq
That is, each $z_n$ is fixed under $\tilde{\alpha}$.
From part (1) of the theorem, we see
$$z_n = (\pi_\tau \otimes h) \alpha(b_n) = \pi_\tau(\bar{b}_n), $$
where
$$\bar{b}_n =  (1 \otimes h) \alpha(b_n) \in B^\alpha$$
is the average of of $b_n$.
Since $\alpha$ is ergodic, $\bar{b}_n$ is a scalar. This implies
that each $z_n$ is also a scalar. Consequently, the operator $z$,
as a limit of the $z_n$'s in the weak operator topology,
is a scalar.
\QED
\vv
{\em Remarks.}
Define on the Hilbert $A$-module $H \otimes A$ (conjugate linear in
the second variable) an operator $u$ by
\bgeq
u(\tilde{b} \otimes a) = (\pi_\tau \otimes 1) \alpha(b)
(\tilde{1}_B \otimes a).
\ndeq
Then one verifies that $u$ is a unitary representation of the quantum group
$A$ (cf \cite{Wor5,BS2,W1}) and $(\pi_\tau, u)$ satisfies the following
covariance condition in the sense of 0.3 of \cite{BS2}:
\bgeq
\label{covariantaction2}
(\pi_\tau \otimes 1) \alpha(b) = u (\pi_\tau(b) \otimes 1) u^*.
\ndeq
The operator $U$ defined above is given by
$$U = (1 \otimes \pi_h) u.$$
The pair $(\pi_\tau, U)$ along with the
relation \eqref{covariantaction1} can be called a covariant system
in the framework of Hopf von Neumann algebras.
Note that part (1) of the above theorem also gives a conceptual
proof of Proposition 4.2.(i) of \cite{Nakagami1}, where a rather complicated
(and nonconceptual) proof is given.
\vv
{\em Notation.} Let $v$ be a unitary representation of
a quantum group $A$ on some finite dimensional Hilbert space $H_v$
\cite{Wor5,BS2,W1}. Define
\bgeq
Ad_v (b) = v ( b \otimes 1) v^*, \; \; b \in B(H_v).
\ndeq
Then using Proposition 3.2 of \cite{Wor5}, we see that $Ad_v$ is an action
of $A$ on $B(H_v)$ (see also the remark after the proof of Theorem 4.1 in
\cite{W15}). It will be called the {\bf adjoint action} of the quantum group
$A$ for the representation $v$. Note that unlike in the case of locally
compact groups, for quantum groups we have in general
\bgeq
Ad_{v \otimes_{in} w} \neq Ad_v \otimes_{in} Ad_w,
\ndeq
where $\otimes_{in}$ denotes the interior tensor product
representations \cite{Wor5,W2}.

For other basic notions on compact quantum groups, we refer the reader
to \cite{Wor5,W1,W2}.

\section{
\label{UHF}
Ergodic actions of $A_u(Q)$ 
on the Powers factor $R_\lambda$ and the
Murray-von Neumann 
factor $R$}

We construct in this section an ergodic action of the universal
quantum group $A_u(Q)$ on the type $\mathrm{III}_\lambda$
Powers factor $R_\lambda$ for a proper choice of $Q$ and an ergodic action
of $A_u(n)$ on the type $\mathrm{II}_1$ Murray-von Neumann factor $R$.
These are obtained as consequences of \thref{main1} below.

Recall \cite{W1,W5,W5'}
that for every non-singular $n \times n$ complex matrix $Q$ ($n > 1$
in the rest of this paper), the universal compact quantum group $(A_u(Q), u)$
is generated by $u_{ij}$ ($i,j = 1 , \cdots, n$) with defining relations
(with $u = ( u_{ij} )$):
\vv
$A_u(Q): \; \; \;
u^* u = I_n = u u^*, \; \; \;
u^t Q {\bar u} Q^{-1} = I_n = Q {\bar u } Q^{-1} u^t$;
\vv
There is also another related family of quantum groups $A_o(Q)$ 
\cite{W1,W5,W5',Banica1}:
\vv
$A_o(Q): \; \; \;
u^t Q u Q^{-1} = I_n = Q u Q^{-1} u^t, \; \; \;
{\bar u} = u$;  (here $Q > 0$)
\vv
Part (1) of the next proposition gives a characterization of $A_u(Q)$
in terms of the functional $\phi_Q$ defined below.

\bgprop
\label{traceQ}
Consider the adjoint action $Ad_u$ corresponding to the fundamental
representation $u$ of the quantum group $(A_u(Q), u)$.

(1). The quantum group $(A_u(Q), u)$ is the
largest compact matrix quantum group such that its
action $Ad_u$ on $M_n( {\Bbb C} )$ leaves invariant the functional
$\phi_Q$ defined by
\bgeqq
\phi_Q (b) = Tr ( Q^t b), \; \; \; b \in M_n({\Bbb C}).
\ndeqq

(2). $Ad_u$ is an ergodic action if and only if
$Q = \lambda E$, where $\lambda$ is a nonzero scalar,
$E$ is the positive matrix $(1 \otimes h) u^t \bar{u}$ (cf \cite{W5}),
 and $h$ is the Haar measure of $A_u(Q)$.
\ndprop
\pf

(1). It is a straightforward calculation to verify that the action
$Ad_u$ of $A_u(Q)$ leaves the functional $\phi_Q$ invariant.

Assume that $(A, v)$ is a compact quantum group such that
$Ad_v$ leaves $\phi_Q$ invariant ($v = (v_{ij})_{i,j=1}^n$). Then the
$v_{ij}$'s satisfy the defining relations for $A_u(Q)$. Hence $(A, v)$ is a
quantum subgroup of $A_u(Q)$.

(2). A matrix $S$ is fixed by $Ad_u$ if and only if $S$ intertwines
the fundamental representation $u$ with itself. Hence the action $Ad_u$
is ergodic if and only if the fundamental representation $u$ is irreducible.
When $Q = \lambda E$, then $A_u(Q) = A_u(E)$.
Since $E$ is positive, $u$ is irreducible (cf. \cite{Banica2}).

On the other hand we have (see \cite{W5})
$$(u^t)^{-1} = E \bar{u} E^{-1},$$
where $E$ is defined as in the proposition.
We also have
$$(u^t)^{-1} = Q \bar{u} Q^{-1}.$$
Hence
$$ E \bar{u} E^{-1} = Q \bar{u} Q^{-1} \; \;
{\mbox and} \; \;
   Q^{-1} E \bar{u} = \bar{u} Q^{-1} E .$$
If $u$ is irreducible, then so is $\bar{u}$ and therefore $Q^{-1} E =$ scalar.
\QED
\vv
{\em Note.} The proof of necessary condition in (2) above was pointed to us
by Woronowicz. Our original proof contains an error.

In general the invariant functional $\phi_Q$ defined above is not
a trace, even if the action $Ad_u$ is ergodic.
However, for ergodic actions of compact groups on operator
algebras, one has the following finiteness theorem of
H$\o$egh-Krohn-Landstad-St$\o$rmer \cite{HLS}:

\bgth
\label{ergodicHLS}
If a von Neumann algebra admits
an ergodic action of a compact group $G$, then

(a). this von Neumann algebra is finite;

(b). the unique $G$-invariant state is
a trace on the von Neumann algebra.
\ndth

The proposition above shows that part (b) of this finiteness
theorem is no longer true for compact quantum groups in general.
We now show that part (a) of the above finiteness theorem is
false for compact quantum groups either: not only can compact quantum
groups act on infinite algebras,
they can act on purely infinite factors (type III factors).

\bgdf
\label{compatiblesystem}
Let $(B_i, \pi_{ji})$ be an inductive system of $C^*$-algebras
($i, j \in I$). For each $i \in I$, let $\alpha_i$ be an action of
a compact quantum group $A$ on $B_i$. We say that the actions
$\alpha_i$ are a {\bf compatible system of actions} for
$(B_i, \pi_{ji})$  if for each pair $i \leq j$, the following
holds,
$$ (\pi_{ji} \otimes 1 ) \alpha_i = \alpha_j \pi_{ji}.$$
\nddf

The following lemmas will be used in the next theorem. Preserve the
notation in \dfref{compatiblesystem}. Let $\pi_i$ be
the natural embedding of $B_i$ into the inductive limit $B$ of the $B_i$'s.

\bglm
Put for each $i \in I$
$$\alpha \pi_i (b_i) = (\pi_i \otimes 1) \alpha_i (b_i),
\; \; \; b_i \in B_i.$$
Then $\alpha$ induces a well defined action of the quantum group $A$ on $B$.
The action $\alpha$ is ergodic if and only if each $\alpha_i$ is.

Assume further that $\phi_i $ is an inductive system of states on
$B_i$ and that each $\phi_i$ is invariant under $\alpha_i$.
Then $\alpha$ leaves invariant the inductive limit state
$\tau = \lim \phi_i$.

\ndlm
\pf
Let $j > i$, so $\pi_{ji}(b_i) \in B_j$.
Then by the formula of $\alpha$ given in the lemma,  we have
\bgeqq
\hspace{2cm}
\alpha \pi_j (\pi_{ji} (b_i) ) = (\pi_j \otimes 1) \alpha_j (\pi_{ji}(b_i)).
\hspace{2cm} (*)
\ndeqq
Since $ \pi_j \pi_{ji} = \pi_i$, the left hand of the above is equal to
$$\alpha \pi_i (b_i) = (\pi_i \otimes 1) \alpha_i (b_i).$$
From the compatibility
condition we see that the right hand side of $(*)$ is equal to
\bgeqq
 (\pi_j \otimes 1) (\pi_{ji} \otimes 1 ) \alpha_i (b_i)
= (\pi_i \otimes 1) \alpha_i (b_i).
\ndeqq
This shows that $\alpha$ is well defined on the dense subalgebra
${\cal B} = \bigcup \pi_i ({\cal B_i})$ of $B$, where ${\cal B}_i$
is the dense *-subalgebra of $B_i$ according to \dfref{qact}.
It also clear that $\alpha$ is bounded and satisfies conditions
of Definition \ref{qact}. Hence $\alpha$ induces a well defined action of
the quantum group on $B$.

Assume that each $\alpha_i$ is ergodic. It is clear that the action
$\alpha$ is ergodic on the dense *-subalgebra $\cal B$.
Hence $\alpha$ is ergodic on $B$ by \propref{denseB}.
%
%

Conversely, if $\alpha$ is ergodic, then the restrictions $\alpha_i$
of $\alpha$ to $B_i$ is clearly ergodic.

We now show that $\tau$ is invariant under $\alpha$.
Note that $\tau (\pi_i(b_i) =  \phi_i (b_i)$. From this we have
\bgeqq
(\tau \otimes 1) \alpha (\pi_i(b_i)) &=&
(\tau \otimes 1) (\pi_i \otimes 1) \alpha_i(b_i)) \\
&=& (\phi_i \otimes 1) \alpha_i(b_i)) = \phi_i (b_i) \\
&=& \tau (\pi_i (b_i) ).
\ndeqq
By density of $\cal B$ in $B$, we have
$$(\tau \otimes 1) \alpha (b) = \tau (b), \; \; \; b \in B.$$
This completes the proof of the lemma.
\QED
\vv
{\em Note.} Not every action of a compact quantum group
on an inductive limit of $C^*$-algebras arises from a compatible
system of actions of $A$.
\bglm
Let $u_k$ be a unitary representation of a compact quantum group $A$
on $V_k$ for each natural number $k$. Assume that $Ad_{u_k}$ leaves
invariant a functional $\psi_k$ on $B(V_k)$. Then
the action $Ad_{u_1 \otimes_{in} \cdots \otimes_{in} u_k}$ leaves
the functional $\phi^k = \psi_1 \otimes \cdots \otimes \psi_k$ invariant.
\ndlm
\pf Straightforward calculation.
\QED
\vv
Let $Q \in GL(n, {\Bbb C})$ be a positive matrix with trace $1$.
We now construct
a sequence of actions $\alpha_k$ of the compact quantum group
$(A_u(Q), u)$ on $M_n({\Bbb C})^{\otimes k}$. Denote by $u^k$ the $k$-th
fold interior tensor product of the representation $u$, i.e.,
$$ u^k = u \otimes_{in} \cdots \otimes_{in} u,$$
see \cite{W2} for the definition of the interior tensor product $\otimes_{in}$.
Put
$$\alpha_k = Ad_{u^k}, \; \; \;
\phi_Q^k = \phi_Q^{\otimes k} = \phi_Q \otimes \cdots \otimes \phi_Q.$$
Let
$$B = \lim_{ k \rightarrow \infty} M_n({\Bbb C})^{\otimes k}, \; \; \;
\tau_Q = \lim_{ k \rightarrow \infty} \phi_Q^k,$$
$$ {\frak B} = \pi_{Q}(B)^{\prime \prime}, \; \; \;
{\frak A} = \pi_h (A_u(Q))^{\prime \prime},
$$
where $\pi_{Q}$ and $\pi_h$ are respectively the GNS-representations for
the positive functional $\tau_Q$ and the Haar state $h$ on $A_u(Q)$.

\bgth
\label{main1}
The actions $\alpha_k$ ($ k = 1, 2, \cdots$) of $A_u(Q)$
forms a compatible system of ergodic actions leaving the functionals
$\phi^k_Q$ invariant.
These actions gives rise to a natural ergodic action on the UHF algebra $B$
leaving invariant the positive functional $\tau_Q$, which in turn
lifts to an ergodic action on the von Neumann algebra $\frak B$.
\ndth
\pf
It is straightforward to verify that the actions $\alpha_k$ are a
compatible system of actions. Since each $u^k$ is irreducible
(cf \cite{Banica2}), we see that the actions $\alpha_k$ are ergodic.
By \ref{traceQ}.(1), $\phi_Q$ is invariant under the action $Ad_u$.
Hence applying the lemmas above, we see that the functionals
$\phi^k_Q$ are invariant under the actions $\alpha_k$, and these actions
gives rise to an ergodic action of the quantum group $A$ on $B$ leaving
$\tau_Q$ invariant.

Now apply \thref{inducedaction},
the action $\alpha$ on $B$ induces an ergodic action
$$\tilde{\alpha}: \; {\frak B} \longrightarrow {\frak B} \otimes {\frak A}$$
at the von Neumann algebra level defined by
$$\tilde{\alpha} (\pi_Q (b)) = (\pi_Q \otimes \pi_h) \alpha(b),$$
where $b \in B$.
\QED

\bgcor
\label{Powers}
Take
$$ Q =
\left(
\begin{array}{cc}
a & 0    \\
0 & 1-a
\end{array}
\right), \; \; a \in (0, 1/2).
 $$ Then $\tau_Q$ is the Powers state, so
 the quantum group $A_u(Q)$ acts ergodically on
 the Powers factor $R_\lambda$ of type $\mathrm{III}_\lambda$,
 where $\lambda = a /(1-a)$.
\ndcor
%
%

\bgcor
\label{Murray-von Neumann}
(compare \cite{AWassermann3})
Take $ Q = I_n$.
Then $\tau_Q$ is the unique trace on the UHF algebra $B$ of type
$n^\infty$, so the quantum group $A_u(n)=A_u(I_n)$ acts ergodically on
 the hyperfinite $\mathrm{II}_1$ factor $R$.
\ndcor

We will see in \secref{Cuntz} that for appropriate choice of $Q$, the
quantum groups $A_u(Q)$ act on the injective factor $R_\infty$ of type
$\mathrm{III}_1$ also.
It would be interesting to know whether compact quantum groups admit
ergodic actions on factors of type $\mathrm{III}_0$ too.

\section{
\label{fixedpoint}
Fixed point subalgebras of quantum subgroups}

In this section, we show that although the actions of the universal quantum
groups $A_u(Q)$ constructed in the last section are ergodic, when restricted
to some of their non-trivial quantum subgroups,
we obtain interesting large fixed point algebras.

Let
$Q =
\left(
\begin{array}{cc}
a & 0    \\
0 & 1-a
\end{array}
\right)$, as in \ref{Powers}.
Put $q = \lambda^{1/2} = (a /(1-a))^{1/2}$.
Then from the definitions of $SU_q(2)$ and $A_u(Q)$,
we see that $SU_q(2)$ is a quantum subgroup of $A_u (Q)$.
By restriction, we obtain from the action of $A_u(Q)$ an action
of $SU_q(2)$ on $R_\lambda$. The fixed point subalgebra of $R_\lambda$
under the action of $SU_q(2)$ is generated by the Jones projections
$\{ 1, e_1, e_2, \cdots, \}$. The restriction of the Powers states $\tau_Q$
to this fixed point algebra is a trace and its values on the
Jones projections gives the Jones polynomial. See the book of
Jones \cite{Jones2}.

Now take $Q = \frac{1}{n} I_n$. We have \corref{Murray-von Neumann}.
For simplicity of notation, let $\tau$ denote the trace $\tau_Q$ on the
UHF algebra $B$. There are two special quantum subgroups of $A_u(n)$:
$SU(n)$ and $A_o(Q) = A_o(n)$.
By 4.7.d. of \cite{GoodmanHarpeJones}, for any closed subgroup $G$ of
$SU(n)$, the fixed point algebra $R^G$ is a II$_1$ subfactor of
$R$. We now show that the same result holds for quantum subgroups
of $A_o(n)$. For this, it suffices to prove the following

\bgprop
The fixed point subalgebra  $R^{A_o(n)}$ of $R$  for the quantum subgroup
$A_o(n)$ of $A_u(n)$ is a $\mathrm{II}_1$ factor and the action
of $A_o(n)$ on $R$ is prime.
\ndprop
\pf
Put $\beta = n^2$. By \cite{Banica1}, the fixed point subalgebra
of $M_n({\Bbb C})^{\otimes k}$ for the action
$$\alpha_k = Ad_{u^k}$$
is generated by $1, e_1, \cdots, e_{k-1}$, where
$u$ is the fundamental representation of the quantum group $A_o(n)$,
$$e_s = I_{H^{\otimes(s-1)}} \otimes \sum_{i,j} \frac{1}{n} e_{ij} \otimes
e_{ij} \otimes I_{H^{\otimes (k-s-1)}},$$
and $H = {\Bbb C}^n$. The $e$'s  satisfies the relations:

(i). $e_s^2 = e_s = e^*_s$;

(ii). $e_s e_t = e_t e_s$, $1 \leq s, t \leq k-1$, $|s-t| \geq 2$;

(iii). $\beta e_s e_t e_s = e_s$, $1 \leq s , t \leq k-1$, $|s-t| = 1$.

We now show that the restriction of $\tau$ on the fixed point subalgebra
of $M_n({\Bbb C})^{\otimes k}$ satisfies the Markov
trace condition of modulus $\beta$, where $\tau$ is the
trace on $R$. Namely, we will verify the identity
         $$\tau ( w e_{k-1} ) = \frac{1}{\beta} \tau (w)$$
for $w$ in the subalgebra of $M_n({\Bbb C})^{\otimes k}$ generated
by $1, e_1, \cdots, e_{k-2}$. By Theorem 4.1.1 and Corollary 2.2.4 of
Jones \cite{Jones1}, this will complete the proof of the proposition.
It will also follow that the action of $A_o(n)$ on $R$ is prime.
To verify this, it suffices by Proposition 2.8.1 of \cite{GoodmanHarpeJones}
to check the Markov trace condition for $w$ of the form
 $$ w = (e_{i_1} e_{i_1-1} \cdots e_{j_1})
        (e_{i_2} e_{i_2-1} \cdots e_{j_2}) \cdots
        (e_{i_p} e_{i_p-1} \cdots e_{j_p}),$$
        where
        \bgeqq
        & & 1 \leq i_1 < i_2 \cdots i_p \leq k-2 \\
        & & 1 \leq j_1 < j_2 \cdots j_p \leq k-2 \\
        & & i_1 \geq j_1, i_2 \geq j_2, \cdots, i_p \geq j_p \\
        & & 0 \leq p \leq k-2
       \ndeqq
If $i_p < k-2$ then it is easy to see that
$$\tau (w e_{k-1}) = \tau (w) \tau(e_{k-1}) = \frac{1}{\beta} \tau (w),$$
noting that $\tau(e_{k-1}) = \frac{1}{\beta}$.
Hence we can assume $i_p = k-2$. Let $l_w$ be the length of the word $w$.
Then $w$ takes the form
$$w = \sum (\frac{1}{n})^{l_w} ( \cdots ) \otimes e_{ab} \otimes 1,$$
where the summation is over the indices $a, b$ and some other indices
that need not to be specified, and the terms in $( \cdots )$ are certain
elements of $M_n({\Bbb C})^{\otimes (k-2)}$ that need not to be
specified either (the components in the tensor product of the terms in
$(\cdots)$ are products of $e_{ij}$'s). We have then
\bgeqq
\tau (w e_{k-1}) &=& (\frac{1}{n})^{l_w + 1}
\sum_{x,y} \sum \tau ((\cdots) \otimes e_{ab} e_{xy} \otimes e_{xy}) \\
 &=& (\frac{1}{n})^{l_w + 1}
\sum_{x,y} \sum \tau ((\cdots) \otimes e_{ab} e_{xy} \otimes 1 )
 \tau (I_{H^{\otimes (k-1)}} \otimes e_{xy}) \\
 &=& (\frac{1}{n})^{l_w + 1}
 \tau (\sum_{x}
 \sum ((\cdots) \otimes e_{ab} e_{xx} \otimes 1 ) \frac{1}{n} ) \\
 &=& \frac{1}{\beta} \tau (w).
 \ndeqq
The proof is complete.
\QED
\vv
{\em Remarks.}
(1).
In view of the above result, fixed point algebras of quantum subgroups
of $A_o(n)$ give examples of subfactors.
Therefore, it would be interesting to classify finite quantum subgroups
of $A_o(n)$ and study them in the light of Jones' theory,
see \cite{GoodmanHarpeJones} for this in the case of the Lie group $SU(2)$.
Note that the quantum $A_o(n)$ contains the quantum permutation
group $A_{aut}(X_n)$ of $n$ point space $X_n$ (see \cite{W15})
and many other interesting quantum subgroups (see \cite{W1}).
It would also be interesting to determine the fixed
point subalgebras of the quantum subgroups of $SU_{-1}(n)$
($SU_{-1}(n)$ is a quantum subgroup of $A_u(n)$
because its antipode has period 2 \cite{W1,W5}).
We refer the reader to Banica \cite{Banica5} for some interesting
related results.

(2).
Note that since
\bgeqq
Ad_{v \otimes_{in} w} \neq Ad_v \otimes_{in} Ad_w,
\ndeqq
for unitary representation $v$ and $w$ of $A_o(n)$,
we do not have a commuting square like the one on p222 of
\cite{GoodmanHarpeJones} for a given quantum subgroup $G$ of $A_o(n)$.

\section{
\label{Cuntz}
Ergodic actions of $A_u(Q)$ on the Cuntz algebra and the injective factor of
type $\mathrm{III}_1$}

%
%

The Cuntz algebra ${\cal O}_n$ is an infinite simple
$C^*$-algebra without trace, hence by \cite{HLS}
it does not admit an ergodic action of a compact group.
Recall that the Cuntz algebra ${\cal O}_n$ is the simple $C^*$-algebra
generated by $n$ isometries $S_k$ ($k=1, \cdots, n$) such that
\bgeq
\sum_{k=1}^n S_k S_k^* = 1.
\ndeq
Just as $U(n)$, the compact matrix quantum group $A_u(Q)$ acts
on ${\cal O}_n$ in a natural manner \cite{DopRob6,Dop1,KMW},
where $Q$ is a positive matrix of trace $1$ in $GL(n, {\Bbb C})$:
\bgeq
\alpha(S_j) = \sum_{i = 1}^n S_i \otimes u_{ij},
\ndeq
the dense *-algebra ${\cal B}$ of \dfref{qact} being
the *-subalgebra $^0{\cal O}_n$ of ${\cal O}_n$
generated by the $S_i$'s, see Doplicher-Roberts \cite{DopRob5}.
However, unlike the actions of compact groups on ${\cal O}_n$, we have
\bgth
\label{main2}
The above action $\alpha$ of the quantum group $A_u(Q)$ on
${\cal O}_n$ is ergodic, the unique $\alpha$-invariant state
on ${\cal O}_n$ is the quasi-free state $\omega_Q$ associated with
$Q$ \cite{Evans1}.
\ndth
\pf
Let $H$ be the Hilbert subspace of ${\cal O}_n$ linearly spanned by
the $S_k$'s. Let $(H^s, H^r)$ be the linear span of elements of the form
$S_{i_1} S_{i_2} \cdots S_{i_r} S^*_{j_s} \cdots S^*_{j_2} S^*_{j_1}$.
Then $ ^0{\cal O}_n $ is the linear span of
all the spaces $(H^s, H^r)$ , $r, s \geq 0$  (see \cite{DopRob5}).
Observe that each of the spaces $(H^s, H^r)$ is invariant
under the action $\alpha$:
$$
\alpha(S_{i_1} S_{i_2} \cdots S_{i_r} S^*_{j_s} \cdots S^*_{j_2} S^*_{j_1})
=  $$
$$
\sum_{k_1, \cdots k_r, l_1, \cdots, l_s = 1 }^n
S_{k_1} S_{k_2} \cdots S_{k_r} S^*_{l_s} \cdots S^*_{l_2} S^*_{l_1}
\otimes u_{k_1 i_1} u_{k_2 i_2} \cdots u_{k_r i_r} u^*_{l_s j_s}
\cdots u^*_{l_2 j_2} u^*_{l_1 j_1}.
$$
Hence $(id \otimes h) \alpha((H^s, H^r))$ is the space of
the fixed elements of $(H^s, H^r)$ under $\alpha$, where $h$ is the
Haar state on $A_u(Q)$. For $r \neq s$, the tensor product representations
$u^{\otimes r}$ and $u^{\otimes s}$ of the fundamental representation
$u$ of the quantum group $A_u(Q)$ are {\em inequivalent and irreducible}
\cite{Banica2}.
Hence by Theorem 5.7 of Woronowicz \cite{Wor5}, for $r \neq s$,
\bgeq
h(u_{k_1 i_1} u_{k_2 i_2} \cdots u_{k_r i_r} u^*_{l_s j_s}
\cdots u^*_{l_2 j_2} u^*_{l_1 j_1}) = 0,
\ndeq
and therefore $(H^s, H^r)$ has no fixed point other than $0$.
For $r = s$, identifying the elements
$$S_{i_1} S_{i_2} \cdots S_{i_r} S^*_{j_r} \cdots S^*_{j_2} S^*_{j_1}$$
of $(H^r, H^r)$ with the matrix units
$$e_{i_1 j_1} \otimes e_{i_2 j_2} \otimes
\cdots \otimes e_{i_r j_r}$$
of $M_n({\Bbb C})^{\otimes r}$, the action $\alpha$ on $(H^r, H^r)$
is identified with the action $\alpha_r$ on
$M_n({\Bbb C})^{\otimes r}$ of \thref{main1}. Hence the fixed elements
of $(H^r, H^r)$ under $\alpha$ are the scalars.
Consequently, the fixed elements of  $^0{\cal O}_n$
under $\alpha$ are the scalars. By
\propref{denseB}, $\alpha$ is ergodic on ${\cal O}_n$.

Let $\phi$ be the (unique) $\alpha$-invariant state
on ${\cal O}_n$. Then for $x \in (H^r, H^s)$ with
$r \neq s$, $r, s \geq 0$, we have
\bgeqq
\phi(x) = h((\phi \otimes 1) \alpha(x))
        = \phi ((1 \otimes h) \alpha(x)).
\ndeqq
But
$ (1 \otimes h) \alpha(x) = 0$
according to the computation above.
Hence $\phi(x) = 0$.
From the consideration of the last paragraph,
$\alpha$ restricts to an ergodic action on the subalgebra
$(H^k, H^k)$ of  ${\cal O}_n$. Identifying $(H^k, H^k)$
with $M_n({\Bbb C})^{\otimes k}$ as above, we see that
$$\phi(x) = \phi_Q^k(x), \; \; \; x \in (H^k, H^k),$$
where $\phi_Q^k$ is the functional in \thref{main1}.
This shows that $\phi$ is the quasi-trace state $\omega_Q$
associated with $Q$ (cf \cite{Evans1}).
\QED
\vv
We can assume that $Q = diag(q_1, q_2, \cdots, q_n)$
is a diagonal positive matrix with trace 1, since
$A_u(Q)$ and $A_u(VQV^{-1})$ are similar to each other \cite{W5}.
Let $\beta$ be a positive number. Define numbers $\omega_1,
\omega_2, \cdots, \omega_n$ by
\bgeq
diag(e^{-\beta \omega_1}, e^{-\beta \omega_2}, e^{-\beta \omega_n})
= diag(q_1, q_2, \cdots, q_n).
\ndeq
Let $\pi_Q$ be the GNS representation of the $\alpha$-invariant
state $\omega_Q$ of ${\cal O}_n$.
Then by Theorem 4.7 of Izumi \cite{Izumi1} and \thref{inducedaction}, we have

\bgcor
If $\omega_1/\omega_k$ is irrational for some $k$, then the compact
quantum group $A_u(Q)$ acts ergodically on the injective
factor $\pi_Q ({\cal O}_n)^{\prime \prime}$ of type $\mathrm{III}_1$.
%
\ndcor
\noindent
{\em Remarks.}
(1).
The big quantum semi-group $U_{nc}(n)$ of Brown also acts on ${\cal O}_n$
in the same way as $A_u(Q)$ on ${\cal O}_n$ above. See Brown \cite{Brown1} and
4.1 of Wang \cite{W1} for the quantum semi-group structure on $U_{nc}(n)$.

(2).
If the $\omega_1/\omega_k$'s are rational for all $k$, then
$\pi_Q ({\cal O}_n)^{\prime \prime}$ is an
injective factor of type $\mathrm{III}_\lambda$,
on which $A_u(Q)$ acts ergodically,
where $\lambda$ is determined from an equation involving $q_1, \cdots, q_n$
(see \cite{Izumi1}). In particular, taking $A_u(Q)=A_u(n)$, we see that even
the compact matrix quantum group $A_u(n)$ of Kac type admits ergodic
actions on both the infinite $C^*$-algebra ${\cal O}_n$ and the
injective factor $\pi_Q({\cal O}_n)^{\prime \prime}$ of type
III$_{\frac{1}{n}}$. In view of \corref{Murray-von Neumann},
it would be interesting to solve the following problem:
\vv
{\bf Problem:}
Does a compact matrix quantum group of non-Kac type admit ergodic action
on the hyperfinite $\mathrm{II}_1$ factor $R$?

\section{
\label{quotient}
Ergodic actions on quotient spaces}

In this section, we study ergodic actions of compact quantum groups on
their quantum quotient spaces. We also give an example to show that,
contrary to the classical situation,
not all ergodic actions arise in this way.

Fix a quantum subgroup $H$ of a compact quantum group $G$, which is given
by a surjective morphism $\theta$ of Woronowicz Hopf $C^*$-algebras
from $C(G)$ to $C(H)$. Let $h_H$ and $h_G$ be respectively the Haar
states on $C(H)$ and $C(G)$. Then there is a natural action $\beta$ of the
quantum group $H$ on $G$ given by
\bgeq
\beta: C(G) \longrightarrow C(H) \otimes C(G), \; \; \;
\beta = (\theta \otimes 1) \Phi_G,
\ndeq
where $\Phi_G$ is the coproduct on $C(G)$. The
quotient space $H \backslash G$ is defined by
the fixed point algebra of $\beta$ (cf \cite{Pod6}):
\bgeq
C(H \backslash G) = C(G)^\beta =
\{ a \in C(G): (\theta \otimes 1) \Phi_G (a)
= 1 \otimes a \} .
\ndeq
The restriction of $\Phi_G$ to $C(H \backslash G)$ defines a natural
action $\alpha$ of $G$ on $C(H \backslash G)$:
\bgeq
\alpha = \Phi_G |_{C(H \backslash G)} : C(H \backslash G)
\longrightarrow C(H \backslash G) \otimes C(G).
\ndeq
The dense *-subalgebras of \dfref{qact} for the actions $\beta$
and $\alpha$ are the natural ones.
Note that $E = (h_H \otimes 1) \beta = (h_H  \theta \otimes 1) \Phi_G $
is a projection of norm
one from $C(G)$ to $C(H \backslash G)$ (cf \propref{denseB} and \cite{Pod6}).

\bgprop
\label{quotients}
In the situation as above, we have

(1). the action $\alpha$ of $G$ on $C( H \backslash G)$ is ergodic;

(2). $C(H \backslash G)$ has a unique $\alpha$
invariant state $\omega$ satisfying
\bgeq
h_G(a) = \omega ((h_H \theta \otimes 1) \Phi_G (a)), \; \; \; a \in C(G).
\ndeq
Namely, $\omega$ is the restriction of $h_G$  on $C(H \backslash G)$.
\ndprop
{\em  Note.} Part (2) of the proposition above is the analogue of the the
following well known integration formula in the classical situation:
 $$\int_G a(g) dg = \int_{H \backslash G} \int_H a(hg) dh d \omega (g) ,
 \; \; \; a \in C(G).$$

\pf
(1). Let $a \in C(H \backslash G)$ be fixed under $\alpha$, i.e.,
$$ \hspace{2cm} \alpha (a) = a \otimes 1. \hspace{2cm} (**)$$
Since $\alpha(a) = \Phi_G (a)$ and since (by the definition of
$C( H \backslash G)$)
$$(\theta \otimes 1) \Phi_G (a) = 1 \otimes a ,$$
it follows that
$$( \theta \otimes 1 ) \alpha (a) = 1 \otimes a.$$
Using $(**)$ for the left hand side of the above, we get
$$\theta(a) \otimes 1 = 1 \otimes a.$$
This is possible only for $a = \lambda \cdot 1$ for some scalar $\lambda$.

(2). The general result of the existence and uniqueness of the
invariant state for an ergodic action is proven in \cite{Boca1}.
For the special situation we consider here, we now not only prove the
existence and uniqueness of the invariant state, but also give the precise
formula of the invariant state.

Let $\omega$ be the restriction of $h_G$ on the subalgebra
$C(H \backslash G)$ of $C(G)$. Since $(h_H \theta \otimes 1) \Phi_G$
is a projection from $C(G)$ onto $C(H \backslash G)$ and $\alpha$ is the
restriction of $\Phi_G$ on $C(H \backslash G)$, the invariance of $\omega$ for
the action $\alpha$ follows from the invariance of the Haar state $h_G$.

Conversely, let $\mu$ be any invariant state on $C(H \backslash G)$.
Using again the fact that $(h_H \theta \otimes 1) \Phi_G$
is a projection from $C(G)$ onto $C(H \backslash G)$,
a standard calculation shows that the functional
$$\phi(a) = \mu ((h_H \theta \otimes 1) \Phi_G (a)), \; \; \; a \in C(G)$$
is a right invariant state, i.e.
$$\phi * \psi (a) = \phi(a), \; \; \; a \in C(G),$$
where $\psi$ is a state on $A$ and
$\phi * \psi = (\phi \otimes \psi) \Phi_G$
is the convolution operation (cf \cite{Wor5}).
From the uniqueness of the Haar state, it follows from this that
$$\phi = h_G, \; \; \; \mu = \omega = h_G|_{C(H \backslash G)}.$$
\QED
\vv
{\em Remarks.}
(1).
Note that the quantum groups $A_u(Q)$, $A_o(Q)$ and $B_u(Q)$ have many
quantum subgroups. In the light of \propref{quotients} and
\thref{inducedaction}, it would be interesting to study
the corresponding operator algebras and the actions on them.
We leave this to a separate work.

(2).
More general than the considerations in \propref{quotients}, if two quantum
groups admit commuting actions on a noncommutative space, then
they act on each other's orbit spaces (not necessarily in an ergodic manner),
just as in the classical situation.
Note that the notion of orbit space corresponds
to fixed point algebra in the noncommutative situation.
\vv
{\em An Example.}
Every transitive action of a compact group $G$ on a
topological space $X$ is isomorphic to the natural action of $G$ on
$H \backslash G$, where $H$ is the closed subgroup of $G$ that fixes
some point of $X$. However, this is no longer true for quantum groups,
even if the space on which the quantum group acts is a classical one.

To see this, let $X_n = \{ x_1, \cdots, x_n \}$ be the space with $n$
points. By Theorem 3.1 of \cite{W15}, the quantum automorphism group
$A_{aut}(X_4)$ of $X_4$ contains the ordinary permutation group $S_4$,
hence it acts ergodically on $X_4$. The quantum subgroup of $A_{aut}(X_4)$
that fix a point, say $x_1$, is isomorphic to $A_{aut}(X_3)$, which is the
same as $C(S_3)$, a (commutative) algebra of dimension $6$.
From \cite{W15}, we know that as a $C^*$-algebra, $A_{aut}(X_n)$
is the same as $C(S_n)$ for $n \leq 3$ and it has
$C^*({\Bbb Z} / 2{\Bbb Z} * {\Bbb Z} / 2{\Bbb Z})$ as a
quotient for $n \geq 4$, where ${\Bbb Z} / 2{\Bbb Z} * {\Bbb Z} / 2{\Bbb Z}$
is the free product of the two-element group
${\Bbb Z} / 2{\Bbb Z}$ with itself, because the entries of the matrix
\bgeqq
\left(
\begin{array}{cccc}
p   & 1-p & 0   & 0    \\
1-p & p   & 0   & 0   \\
0   & 0   & q   & 1-q  \\
0   & 0   & 1-q & q
\end{array}
\right)
\ndeqq
satisfy the commutation relations of the algebra $A_{aut}(X_4)$,
where $p, q$ are the projections generating the $C^*$-algebra
$C^*({\Bbb Z} / 2{\Bbb Z} * {\Bbb Z} / 2{\Bbb Z})$: $p=(1-u)/2$ and
$q=(1-v)/2$, $u$ and $v$ being the unitary generators of the
first and second copies of ${\Bbb Z}/2{\Bbb Z}$ in the free product
(cf \cite{RaeburnSinclair}).

For simplicity of notation, let $C(G) = A_{aut}(X_4)$,
and let ${\cal M}$ be the canonical dense subalgebra of $C(G)$
generated by the coefficients of the fundamental representation
of $G$ (see \cite{Wor5}). Let $H=S_3$, the subgroup of $G$ that
fixes $x_1$, and let ${\cal H} = C(H)$. Let
$\theta$ be the surjection from $C(G)$ to $C(H)$ that embeds
$H$ as subgroup of $G$ (cf. \cite{W15}).
Let $\beta$ be the action defined in the beginning of this section.
We claim that the coset space $H \backslash G$ is not isomorphic to $X_4$ as
a $G$-space (see Sect. 2 of \cite{W15} for the notion of morphism).
Namely, we have

\bgprop
The $G$-algebras $C(H \backslash G)$
(which is defined to be $C(G)^\beta$) and $C(X_4)$ are not isomorphic to
each other.
\ndprop
\pf
Since $C(X_4)$ has dimension $4$,
it suffices to show that  $C(H \backslash G)$ is infinite dimensional.

We make ${\cal M}$ into a
Hopf ${\cal H}$-module (i.e. a compatible system of a left ${\cal H}$
comodule and a left ${\cal H}$ module) as follows.
The restriction of $\beta$ to ${\cal M}$ clearly defines a
left ${\cal H}$ comodule structure:
\bgeq
\beta: {\cal M} \longrightarrow {\cal H} \otimes {\cal M}.
\ndeq
The left ${\cal H}$ module structure on ${\cal M}$ is the trivial one
defined by
\bgeq
& & {\cal H} \otimes {\cal M} \longrightarrow {\cal M}, \\
& & h \cdot m = \epsilon(h)m, \; \; \; h \in {\cal H}, \;  m \in {\cal M}.
\ndeq
By Theorem 4.1.1 of Sweedler \cite{Sweedler}, we have
an isomorphism of left ${\cal H}$ modules
\bgeq
& & {\cal H} \otimes {\cal M}^\beta \cong {\cal M}.
\ndeq
That is
\bgeq
& & {\cal H} \otimes {\cal A}(H \backslash G) \cong {\cal M} \\
& & h \otimes m' \mapsto h \cdot m',
\; \; \; h \in {\cal H}, \; m' \in {\cal A}(H \backslash G),
\ndeq
where ${\cal A}(H \backslash G)= {\cal M}^\beta$ is the canonical dense
subalgebra of $C(H \backslash G)=C(G)^\beta$. Since ${\cal M}$ is infinite
dimensional and ${\cal H}$ is finite dimensional, ${\cal A}(H \backslash G)$
and therefore $C(H \backslash G)$ are also infinite dimensional.
\QED
\vv
{\bf Acknowledgement.}
The author is indebted to Marc A. Rieffel for continual support.
Part of this paper was written while the author was a
member at the IHES during the year July, 1995-Aug, 1996. He
thanks the IHES for its financial support and hospitality during this period.
The author also wishes to thank the Department of
Mathematics at UC-Berkeley for its support and hospitality
while the author holds an NSF Postdoctoral Fellowship there
during the final stage of this paper.

\vspace{1.5cm}
\hfill 
revised Dec 10, 1998


\begin{thebibliography}{99}


\bibitem{BS2} Baaj, S. and Skandalis, G.:
{\rm Unitaires multiplicatifs et dualit\'e pour les produits
crois\'es de $C^*$-alg\`ebres,}
{\em Ann. Sci. Ec. Norm. Sup.} {\bf 26} (1993), 425-488.

\bibitem{Banica1} Banica, T.:
{\rm Th\'eorie des repr\'esentations du groupe quantique
compact libre $O(n)$,}
{\em C. R. Acad. Sci. Paris} t. {\bf 322}, Serie I (1996), 241-244.

\bibitem{Banica2} Banica, T.:
{\rm Le groupe quantique compact libre $U(n)$,}
{\em Commun. Math. Phys.} {\bf 190} (1997), 143-172.

\bibitem{Banica5} Banica, T.:
Quantum groups acting on $n$ points, complex Hadamard matrix and a
construction of subfactors,
math/9806054

%
\bibitem{Boca1} Boca, F.:
{\rm Ergodic actions of compact matrix pseudogroups on $C^*$-algebras,}
in {\em Recent Advances in Operator Algebras,}
{\em Ast\'erisque} {\bf 232} (1995), 93-109.

\bibitem{Brown1} Brown, L.:
{\rm Ext of certain free product of $C^*$-algebras,}
{\em J. Operator Theory} {\bf 6} (1981), 135-141.

\bibitem{CDPR1} Ceccherini, T., Doplicher, S., Pinzari, C. and Roberts, J.E.:
{\rm A generalization of the Cuntz algebras and model actions},
{\em J. Funct. Anal.} {\bf 125} (1994), 416-437.

%
\bibitem{Cuntz} Cuntz, Joachim:
{\rm Simple $C\sp*$-algebras generated by isometries,}
{\em Comm. Math. Phys.} {\bf 57} (1977), no.~2, 173--185.

\bibitem{Cuntz2} Cuntz, Joachim:
{\rm Regular actions of Hopf algebras on the $C\sp \ast$-algebra
generated by a Hilbert space,}
in {\it Operator algebras, mathematical physics,
and low-dimensional topology(Istanbul, 1991)}, 87--100,
A K Peters, Wellesley, MA, 1993; MR 94m:461

\bibitem{Dop1} Doplicher, S.:
{\rm Abstract compact group duals, operator algebras and quantum
field theory,}
{\em Proc. ICM}-1990, Kyoto, Springer, 1991.

\bibitem{DopRob5} Doplicher, S. and Roberts, J.E.:
{\rm Duals of compact Lie groups realized in the Cuntz algebras and
their actions on $C^*$-algebras,}
{\em J. Funct. Anal.} {\bf 74} (1987), 96-120.

\bibitem{DopRob6} Doplicher, S. and Roberts, J.E.:
{\rm Compact group actions on $C^*$-algebras,}
{\em J. Operator Theory} {\bf 19} (1988), 283-305.

%
\bibitem{Evans1} Evans, David E.:
{\rm On ${\cal O}_n$},
{\em Publ. RIMS. Kyoto Univ.} {\bf 16} (1980), 915-927.

\bibitem{GoodmanHarpeJones}
Goodman, F.M. and de la Harpe, P. and Jones, V.F.R.:
{\em Coxeter Graphs and Towers of Algebras,}
MSRI Publ. 14, Springer-Verlag, 1989.

\bibitem{HLS} H$\o$egh-Krohn, R. and Lanstad, M.B. and St$\o$rmer, E.:
{\rm Compact ergodic groups of automorphisms,}
{\em Ann. of Math.} {\bf 114} (1981), 75-86

\bibitem{Izumi1} Izumi, Masaki:
{\rm Subalgebras of infinite $C\sp *$-algebras with finite
Watatani indices. I. Cuntz algebras,}
{\em Comm. Math. Phys.} {\bf 155} (1993),
no.~1, 157--182; MR 94e:46104

\bibitem{Jones1} Jones, V. F. R.:
{\rm Index for Subfactors,}
{\em Invent. Math.} {\bf 72} (1983), 1-25.

\bibitem{Jones2} Jones, V. F. R.:
{\em Subfactors and Knots,}
Regional Conference Series {\bf 80},
Amer. Math. Soc., 1991.

\bibitem{KMW} Konishi, Y., Nagisa, M. and Watatani, Y.:
{\rm Some remarks on actions of compact matrix quantum groups on
$C^*$-algebras,}
{\em Pacific J. Math.} {\bf 153} (1992), 119-127.

\bibitem{Marciniak1} Marciniak, M.:
{\rm Actions of compact quantum groups on $C^*$-algebras},
{\em Proc. AMS} {\bf 126} (1998), 607-616.

%
\bibitem{Nakagami1} Nakagami, Y.:
{\rm Takesaki duality for the crossed product by quantum groups,}
in {\em Quantum and Non-Commutative Analysis,} H. Araki ed.,
Kluwer Academic Publishers, 1993, 263-281.

%
%
\bibitem{Paolucci1} Paolucci, A.:
{\rm Coactions of Hopf algebras on Cuntz algebras and
their fixed point algebras},
{\em Proc. AMS} {\bf 125} (1997), 1033-1042.

%
\bibitem{Pod6} Podles, P.:
{\rm Symmetries of quantum spaces. Subgroups and
quotient spaces of quantum $SU(2)$ and $SO(3)$ groups,}
{\em Commun. Math. Phys.} {\bf 170} (1995), 1-20.

\bibitem{RaeburnSinclair} Raeburn, I. and Sinclair, A.M.:
{\rm The $C^*$-algebra generated by two projections},
{\em Math. Scand.} {\bf 65} (1989), 278-290.


\bibitem{Sweedler} Sweedler, M.E.:
{\em Hopf Algebras,}
Benjamin, New York,  1969.

\bibitem{Daele4} Van Daele, A.:
{\rm The Haar measure on a compact quantum group,}
{\em Proc. Amer. Math. Soc.} {\bf 123} (1995), 3125-3128.

\bibitem{W5} Van Daele, A. and Wang, S. Z.:
{\rm Universal quantum groups,}
{\em International J. Math} {\bf 7}:2 (1996), 255-264.

\bibitem{W1} Wang, S. Z.:
{\rm Free products of compact quantum groups,}
{\em Commun. Math. Phys.} {\bf 167} (1995), 671-692.

\bibitem{W2} Wang, S. Z.:
{\rm Tensor products and crossed products of compact quantum groups,}
{\em Proc. London Math. Soc.} {\bf 71} (1995), 695-720.

\bibitem{W5'} Wang, S. Z.:
{\rm New classes of compact quantum groups,}
Lecture notes for talks at the University of Amsterdam and
the University of Warsaw, January and March, 1995.

\bibitem{W7} Wang, S. Z.:
{\rm Problems in the theory of quantum groups,}
in {\em Quantum Groups and Quantum Spaces},
Banach Center Publication 40 (1997),
Inst. of Math., Polish Acad. Sci.,
Editors: R. Budzynski, W. Pusz, and S. Zakrzewski.
pp67-78

\bibitem{W15} Wang, S. Z.:
{\rm Quantum symmetry groups of finite spaces},
{\em Commun. Math. Phys.} {\bf 195}:1 (1998), 195-211.

\bibitem{AWassermann1} Wassermann, A.:
{\rm Ergodic actions of compact groups on operator algebras I:
General theory,}
{\em Ann. of Math.} {\bf 130} (1989), 273-319.

\bibitem{AWassermann3} Wassermann, A.:
{\rm Ergodic actions of compact groups on operator algebras III:
Classification for $SU(2)$,}
{\em Invent. Math.} {\bf 93} (1988), 309-355.

\bibitem{AWassermann4} Wassermann, A.:
{\rm Coactions and Yang-Baxter equations for ergodic actions and subfactors,}
in {\em Operator Algebras and Applications,} no 2, ed. by D. Eavans and
M. Takesaki, London Math. Soc. Lecture Notes {\bf 136} (1988), 203-236.

\bibitem{Wor4} Woronowicz, S. L.:
Twisted $SU(2)$ group. An example of noncommutative differential calculus,
{\em Publ. RIMS, Kyoto Univ.} {\bf 23} (1987), 117-181.

\bibitem{Wor5} Woronowicz, S. L.:
Compact matrix pseudogroups,
{\em Commun. Math. Phys.} {\bf 111} (1987), 613-665.

\bibitem{Wor6} Woronowicz, S. L.:
Tannaka-Krein duality for compact matrix pseudogroups. Twisted $SU(N)$ groups,
{\em Invent. Math.} {\bf 93} (1988), 35-76.


\end{thebibliography}
\end{document}